\documentclass[12pt]{article}
\usepackage{amssymb}
\usepackage{amsmath}
\usepackage{latexsym}

\newcommand {\R}	{{\textbf R}}
\newcommand {\Mat}	{\text{Mat}}
\def\msn{\medskip\noindent}

\newtheorem{theorem}          {Theorem}

\newtheorem{algorithm}      [theorem]  {Algorithm}

\begin{document}
\title{Universal algorithms for generalized 
discrete matrix Bellman 
equations with symmetric Toeplitz matrix\thanks{%
Research supported by RFBR grant 05-01-00824 and
joint RFBR/CNRS grant 05-01-02807.}}
\author{Serge\u{\i} Sergeev\thanks%
{Department of Physics, Sub-Department of Quantum Statistics and Field Theory,
M.V. Lomonosov Moscow State University, 119992 Leninskie Gory, GSP-2, Moscow,
Russia. E-mail: sergiej@gmail.com}}
\maketitle  

\begin{abstract}
\msn This paper presents two universal algorithms for generalized discrete
matrix Bellman equations with symmetric Toeplitz matrix. The algorithms are
semiring extensions of two well-known methods solving
Toeplitz systems in the ordinary linear algebra.
\end{abstract}

\section{Introduction}
\label{s:intro}
\msn As observed
by B.A.~Carr\'{e}\cite{Car-71,BC-75}, the Gaussian
elimination without pivoting can be viewed as a prototype
for some algorithms on graphs. 
M.~Gondran \cite{Gon-75} and G.~Rote \cite{Rot-85} made this
observation precise by proving 
that the Gaussian elimination, under certain conditions, can
be applied to the linear systems of equations over semirings. 

\msn The notion of universal algorithm over semiring was introduced 
by G.L.~Litvinov, V.P.~Maslov and E.V.~Maslova in \cite{LM-98,LMa-00<Univ>}. 
These papers are to be considered in the framework of
publications
\cite{Mas:87<MO>,KM:97<Id>, 
LMSz-01<IFA>,LSz-05<Tambov>,LSob-00<DAN>,LSob-01<IIA&OP>}) of
the Russian idempotent school, and more generally, in the framework of
idempotent and tropical mathematics, see \cite{LM:05,BCOQ:92,C-G:79<MA>}
and references therein.  
Essentially, an algorithm is called universal
if it does not depend on the computer representation of data and
on a specific realization of algebraic operations involved
in the algorithm \cite{LMa-00<Univ>}.
Linear algebraic universal algorithms include generalized 
bordering method, LU- and
LDM-decompositions for solving matrix equations. These methods are
basically due to B.A.~Carr\'{e}, see also \cite{LMa-00<Univ>}.

\msn It was observed in \cite{LM-98,LMa-00<Univ>} that 
universal algorithms can be implemented by means of
objective-oriented programming supported by
C++, MATLAB, Scilab, Maple and other computer systems and languages.  
Such universal programs can be instrumental in many areas including
the problems of linear algebra, optimization theory, and interval
analysis over positive semirings, see 
\cite{LM-98,LSob-00<DAN>,LSob-01<IIA&OP>,LP-05}.

\msn This paper presents new universal algorithms based on the methods
of Durbin and Levinson, see \cite{GvL:99}, Sect. 4.7. These algorithms
solve systems of linear equations with symmetric Toeplitz matrices.
Our universal algorithms have the same computational complexity $O(n^2)$
as their prototypes which beats
the complexity $O(n^3)$ of the LDM-decomposition method. All algorithms
are described as MATLAB-programs, meaning that they can be actually 
implemented.

\msn The author is grateful to G.L.~Litvinov and A.N.~Sobolevski\u{\i}
for drawing his attention to this problem and for valuable discussions.
  
\section{Semirings and universal algorithms}
\label{s:semirings}
\msn A set $S$ equipped with {\em addition} $\oplus$ and 
{\em multiplication} $\odot$
is a {\em semiring (with unity)} if the following axioms hold:
\begin{itemize}
\item[1)] $(S,\oplus)$ is a commutative semigroup with neutral element 
$\textbf{0}$;
\item[2)] $(S,\odot)$ is a semigroup with neutral element
${\bf 1}\ne {\bf 0}$;
\item[3)] $a\odot(b\oplus c)=(a\odot
b)\oplus(a\odot c)$, $(a\oplus b)\odot c=(a\odot c)\oplus (b\odot c)$
for all $a,b,c\in S$ (distributivity);
\item[4)] ${\bf
0}\odot a=a\odot{\bf 0}={\bf 0}$ for all
$a\in S$.
\end{itemize}
In the sequel, we omit the notation $\odot$ whenever this is convenient
and does not lead to confusion. 

\msn The semiring $S$ is called {\em idempotent} if $a\oplus a=a$ for 
any $a\in S$. In this case $\oplus$ induces the {\em canonical partial
order} relation
\begin{equation}
\label{e:idord}
a\preceq b\Leftrightarrow a\oplus b=b.
\end{equation}

\msn The semiring $S$ is called {\em complete} (cf. \cite{Gol:00}), 
if any subset $\{x_{\alpha}\}\subset S$ is summable
and the infinite distributivity
\begin{equation}
\label{e:infdistr}
\begin{array}{l@{{}={}}l}
c\odot(\bigoplus_{\alpha}x_{\alpha}) &
\bigoplus_{\alpha}(c\odot x_{\alpha}),\\
(\bigoplus_{\alpha}x_{\alpha})\odot c &
\bigoplus_{\alpha}(x_{\alpha}\odot c).
\end{array}
\end{equation} 
holds for all $c\in S$ and $\{x_{\alpha}\}\subset S$.
This property is natural in idempotent semirings and also 
in the theory of partially ordered spaces (cf. G.~Birkhoff 
\cite{Bir:84<LT>}) with partial order (\ref{e:idord}).
Complete idempotent semirings are called 
$a$-complete (cf. \cite{LMSz-01<IFA>}).

\msn Consider the   
{\em closure operation}
\begin{equation}
\label{e:clseries}
a^*=\bigoplus_{i=0}^{\infty} a^i.
\end{equation}
In the complete semirings it is defined for all elements.
The property
\begin{equation}
\label{e:quasiinv}
a^*=\textbf{1}\oplus aa^*=\textbf{1}\oplus a^*a,
\end{equation}
reveals that the closure operation is a natural extension
of $(1-a)^{-1}$.

\msn We give some examples of semirings living on the set
of reals $\R$ totally ordered by $\leq$: the semiring $\R_+$
with customary operations $\oplus=+$, $\odot=\cdot$ and neutral
elements
$\textbf{0}=0$ and $\textbf{1}=1$; the semiring
$\R_{\max}=\R\cup\{-\infty\}$ with operations
$\oplus=\max$ и $\odot=+$, and neutral elements $\textbf{0}=-\infty$,
$\textbf{1}=0$; the semiring $\widehat{\R}_{\max}=\R_{\max}\cup\{\infty\}$,
which is a completion of $\R_{\max}$ with the element $\infty$ satisfying
$a\oplus\infty=\infty$ for all $a$, 
$a\odot\infty=\infty\odot a=\infty$ for $a\ne\textbf{0}$ and 
$\textbf{0}\odot\infty=\infty\odot\textbf{0}=\textbf{0}$; the semiring
$\R_{\max,\min}=\R\cup\{\infty\}\cup\{-\infty\}$ with $\oplus=\max$,
$\odot=\min$, $\textbf{0}=-\infty$, and 
$\textbf{1}=\infty$. 

\msn Consider operation (\ref{e:clseries})
for the examples above. In $\R_+$ the closure $a^*$ equals  
$(1-a)^{-1}$ if $a<1$ and is undefined otherwise;
in $\R_{\max}$ it equals
$\textbf 1$ if $a\leq\textbf{1}$ and is undefined otherwise;
in $\widehat{\R}_{\max}$ we have $a^*=\textbf{1}$
for $a\leq\textbf{1}$ and $a^*=\infty$ for $a>\textbf{1}$; 
in $\R_{\max,\min}$
we have $a^*=\textbf{1}$ for all $a$.
Note that $\widehat{\R}_{\max}$ and $\R_{\max}$ are $a$-complete, so
the closure is defined for any element of these semirings.

\msn The {\em matrix} operations $\oplus$ and $\odot$ 
are defined analogously to their counterparts in linear algebra.
Denote by $\Mat_{mn}(S)$ the set of all $m\times n$ matrices over
the semiring $S$. By $I_n$ we denote the $n\times n$ unity matrix,
that is, the matrix with $\textbf 1$ on the diagonal and
$\textbf 0$ off the diagonal. As usual, we have
$AI_n=I_nA=A$ and $A^0=I_n$ for any
$A\in\Mat_{nn}(S)$. The set $\Mat_{nn}(S)$ of all
$n\times n$ square matrices is a
semiring. Its unity is $I_n$ and its zero is $\textbf{0}_n$, the
square matrix with all entries equal to $\textbf 0$.
If $S$ is complete and/or idempotent, then so is the semiring
$\Mat_{nn}(S)$.  
If $S$ (and hence $\Mat_{nn}(S)$)
is complete, the closure $A^*$ is defined for
any matrix $A$ and it satisfies (\ref{e:quasiinv}). 
Note that if $S$ is partially ordered, then
$\Mat_{mn}(S)$ is ordered elementwise: $A\preceq B$ iff 
$A_{ij}\preceq B_{ij}$ for all $i=1,\ldots,m$ and $j=1,\ldots,n$. If
$S$ is idempotent and canonically ordered (\ref{e:idord}), then 
the elementwise order of $\Mat_{mn}(S)$ also satisfies 
(\ref{e:idord}).

\msn The closure operation of matrices is important for
the {\em (discrete stationary) matrix Bellman equations}
\begin{equation}
\label{e:bellman}
X=AX\oplus B.
\end{equation}
If the closure of $A$  
exists and (\ref{e:quasiinv}) holds, then $X=A^*B$
is a solution to (\ref{e:bellman}). In $a$-complete
idempotent semirings the matrix $A^*B$
is the least solution of this equation with respect to
(\ref{e:idord}). 

\msn Since $A^*$ is a generalization of $(I-A)^{-1}$, the known universal
algorithms for $A^*$ are generalizations of the methods for 
matrix inverses, and the known algorithms for Bellman equations are
generalizations of the methods for $AX=B$. 
Further we consider the generalized bordering 
method.
  
\msn Let $A$ be a square matrix. Closures of its main submatrices $A_k$ can
be found inductively. The base of induction is $A_1^*$, the closure of the
the first diagonal entry. Generally, we represent $A_{k+1}$ as
\[
A_{k+1}=
\left(
\begin{array}{ll}
A_k & g_k\\
h_k^T & a_{k+1}
\end{array}
\right),
\]
assuming that we have found the closure of $A_k$.
In this representation, $g_k$ and $h_k$ are columns with $k$ entries
and $a_{k+1}$ is a scalar. We also represent $A_{k+1}^*$ as
\[
A_{k+1}^*=
\left(
\begin{array}{ll} 
U_k & v_k\\ 
w_k^T & u_{k+1}
\end{array}
\right).
\] 
Using (\ref{e:quasiinv}) we obtain that
\begin{equation}
\label{e:bordrec}
\begin{array}{r@{{}={}}l}
u_{k+1} & (h_k^T A_k^* g_k\oplus a_{k+1})^*,\\ 
v_k & A_k^*g_k u_{k+1},\\
w_k^T & u_{k+1} h_k^T A_k^*,\\
U_k & A_k^* g_k u_{k+1} h_k^T A_k^*\oplus A_k^*.
\end{array}
\end{equation}

\msn Consider the bordering method for finding
the solution 
$x=A^*b$ to (\ref{e:bellman}),
where $X=x$ and $B=b$ are column vectors. Firstly, we have $x^{(1)}=
A_1^* b_1$. Let $x^{(k)}$ be the vector found after $(k-1)$ steps,
and let us write
\[
x^{(k+1)}=
\left(
\begin{array}{l} 
z\\ 
x_{k+1} 
\end{array}
\right).
\]
Using~(\ref{e:bordrec}) we obtain that
\begin{equation}
\label{e:bbordrec}
\begin{array}{r@{\;}c@{\;}l}
x_{k+1} &=& u_{k+1}(h_k^T x^{(k)}\oplus b_{k+1}),\\
z &=& x^{(k)}\oplus A_k^* g_k x_{k+1}^{(k+1)}.
\end{array}
\end{equation}

\msn We have to compute
$A_k^* g_k$. In general it makes a problem, 
but not in the case of the next section when $A$ is symmetrical
Toeplitz.  

\msn We also note that the bordering method described by
(\ref{e:bordrec}) and (\ref{e:bbordrec}) is valid more generally
over Conway semirings, see \cite{Gol:00} for the definition.

\section{Universal algorithms for Toeplitz linear systems} 
\label{s:toeplitz}
\msn Formally, a matrix $A\in\Mat_{nn}(S)$
is called {\em (generalized) Toeplitz} 
if there exist scalars
$r_{-n+1},\ldots ,r_0,\ldots ,r_{n-1}$ such that $A_{ij}=r_{j-i}$
for all $i$ and $j$.  
Informally, Toeplitz matrices are such that their entries are 
constant along any line
parallel to the main diagonal (and along the main diagonal itself).
For example,
\begin{equation}
A=
\left(
\begin{array}{llll}
r_0 & r_1 & r_2 & r_4\\
r_{-1} & r_0 & r_1 & r_2\\
r_{-2} & r_{-1} & r_0 & r_1\\
r_{-3} & r_{-2} & r_{-1} & r_0
\end{array}
\right)
\end{equation}
is Toeplitz.
Such matrices are not necessarily symmetric. However, 
they are always {\em persymmetric}, that is, symmetric with respect
to the inverse diagonal. This property is algebraically expressed as
$A=E_n A^T E_n$, where
$E_n=[e_n,\ldots ,e_1]$. By $e_i$ we denote  
the column whose $i$th entry is $\textbf{1}$ and other entries 
are $\textbf{0}$. The property $E_n^2=I_n$ (where $I_n$ is the
$n\times n$ identity matrix) implies that the product of two
persymmetric matrices is persymmetric. Hence any degree of a persymmetric
matrix is persymmetric, and so is the closure of a persymmetric
matrix. Thus, if $A$ is persymmetric, then 
\begin{equation}
\label{pers}
E_n A^*=(A^*)^T E_n.
\end{equation}

\msn Further we deal only with symmetric Toeplitz
matrices. Consider the equation  
$y=T_n y\oplus r^{(n)}$, where $r^{(n)}=(r_1,\ldots r_n)^T$
and $T_n$ is defined by the scalars $r_0,r_1,\ldots ,r_{n-1}$ so that
$T_{ij}=r_{\mid j-i\mid}$ for all $i$ and $j$. This is a generalization
of the Yule-Walker problem \cite{GvL:99}.  
Assume that we have obtained a
solution $y^{(k)}$ to the system $y=T_k y\oplus r^{(k)}$ 
for some $k$ such that
$1\leq k\leq n-1$, where  $T_k$ is the main $k\times k$ 
submatrix of $T_n$. We write $T_{k+1}$ as
\[
T_{(k+1)}
=
\left(
\begin{array}{ll}
T_k & E_k r^{(k)}\\
r^{(k)T} E_k & r_0
\end{array}
\right).
\]
We also write $y^{(k+1)}$ and $r^{(k+1)}$ as
\[
\begin{array}{ll}
y^{(k+1)}=\left(\begin{array}{l} z\\ \alpha_k\end{array}\right), &
r^{(k+1)}=\left(\begin{array}{l} r^{(k)}\\ r_{k+1}\end{array}\right).
\end{array}
\]
Using (\ref{e:bbordrec}), (\ref{pers}) and
the identity $T_k^* r^{(k)}=y^{(k)}$,
we obtain that
\begin{equation}
\label{e:badjuwrec}
\begin{array}{r@{\;}c@{\;}l}
\alpha_k &=& (r_0\oplus r^{(k)T} y^{(k)})^*(r^{(k)T}E_k y^{(k)}\oplus r_{k+1}),\\
z &=& E_k y^{(k)}\alpha_k\oplus y^{(k)}.
\end{array}
\end{equation}
Denote $\beta_k=r_0\oplus r^{(k)T} y^{(k)}$. The following
argument shows that $\beta_k$
can be found recursively if
$(\beta_{k-1}^*)^{-1}$ exists.
\begin{equation}
\label{e:arg}
\begin{array}{r@{{}={}}l}
\beta_k & r_0\oplus [r^{(k-1)T}\quad r_k]
\left(\begin{array}{l}
E_{k-1} y^{(k-1)}\alpha_{k-1}\oplus y^{(k-1)}\\
\alpha_{k-1}\end{array}\right)=\\
& r_0\oplus r^{(k-1)T} y^{(k-1)}\oplus 
(r^{(k-1)T} E_{k-1} y^{(k-1)}\oplus r_k)\alpha_{k-1}=\\
& \beta_{k-1}\oplus (\beta_{k-1}^*)^{-1}\odot(\alpha_{k-1})^2.
\end{array}
\end{equation}

The argument above is not always valid and this will make us
write two versions of our algorithm, the first one involving
(\ref{e:arg}) and the second one not involving it. We will write these
two versions in one program and 
mark the expressions which refer only to the first
or only to the second  
version by the MATLAB-style comments $\%1$ or $\%2$, respectively.
Collecting the expressions for
$\beta_k$,$\alpha_k$ and $z$, we obtain the
following recursive expression for $y^{(k)}$:
\begin{equation}
\label{e:durbin}
\begin{array}{r@{{}={}}l}
\beta_k & r_0\oplus r^{(k)T} y^{(k)},\qquad \% 2\\
\beta_k & \beta_{k-1}\oplus (\beta_{k-1}^*)^{-1}\odot(\alpha_{k-1})^2,\qquad \% 1\\
\alpha_k & (\beta_k)^*\odot ((r^{(k)T}E_k y^{(k)}\oplus r_{k+1}),\\
y^{(k+1)} & \left(\begin{array}{l}
E_{k} y^{(k)}\alpha_{k}\oplus y^{(k)}\\
\alpha_{k}\end{array}\right).
\end{array}
\end{equation}

Recursive expression (\ref{e:durbin}) is a generalized version  
of the Durbin method for the Yule-Walker problem
\cite{GvL:99}. Using this expression we obtain the following algorithm.
\begin{algorithm}
\label{a:juw}
The Yule-Walker problem for the Bellman equations 
with symmetric Toeplitz matrix.
\end{algorithm}
\noindent
${\bf function}\ y={\bf durbin}(r_0,r)\\
n={\bf size}(r)+1\\
y(1)=r_0^* \odot r(1)\\
\beta=r_0\qquad \% 1\\
\alpha=r_0^*\odot r(1)\qquad \% 1\\
{\bf for}\quad k=1:n-1\\
\-\beta=r_0\oplus r(1:k)\odot y(1:k)\qquad \% 2\\
\-\beta=\beta\oplus (\beta^*)^{-1}\odot\alpha^2\qquad \% 1\\
\-\alpha=\beta^*\odot(r(k:-1:1)\odot y(1:k)\oplus r(k+1))\\
\- z(1:k)=y(1:k)\oplus\alpha\odot y(k:-1:1)\\
\- y(1:k)=z(1:k)\\
\- y(k+1)=\alpha\\
{\bf end}$

Now we consider the problem of finding
$x^{(n)}=T_n^* b^{(n)}$ where $T_n$ is as above
and $b^{(n)}=(b_1,\ldots ,b_n)^T$ is arbitrary.
We also introduce the column $y^{(n)}$ which solves
the Yule-Walker problem:
$y^{(n)}=T_n^* r^{(n)}$. The main idea is to find the expression for
$x^{(k+1)}=T_{k+1}^* b^{(k+1)}$
involving $x^{(k)}$ and $y^{(k)}$.
We write $x^{(k+1)}$ and $b^{(k+1)}$ as
\[
\begin{array}{ll}
x^{(k+1)}=\left(\begin{array}{l} v\\ \mu_k\end{array}\right), &
b^{(k+1)}=\left(\begin{array}{l} b^{(k)}\\ b_{k+1}\end{array}\right).
\end{array}
\]

Making use of the persymmetry of $T_k^*$ and of the identities
$T_k^* b_k=x^{(k)}$ and $T_k^* r_k=y^{(k)}$, 
we specialize expressions~(\ref{e:bbordrec}) and obtain that
\begin{equation}
\label{e:mukv}
\begin{array}{r@{{}={}}l}
\mu_k &(r_0\oplus r^{(k)T} y^{(k)})^*\odot ((r^{(k)T}E_k x^{(k)}\oplus b_{k+1}),\\
v & E_{k} y^{(k)}\mu_{k}\oplus x^{(k)}.
\end{array}
\end{equation}

The coefficient $r_0\oplus r^{(k)T} y^{(k)}=\beta_k$ is again to be 
expressed as $\beta_k=\beta_{k-1}\oplus
(\beta_{k-1}^*)^{-1}\odot(\alpha_{k-1})^2$, if the
closure $(\beta_{k-1})^*$ is invertible.
Using this we obtain the following recursive expression:
\begin{equation}
\label{e:levinson}
\begin{array}{r@{{}={}}l}
\beta_k & r_0\oplus r^{(k)T} y^{(k)},\qquad \% 2\\
\beta_k & \beta_{k-1}\oplus (\beta_{k-1}^*)^{-1}\odot(\alpha_{k-1})^2,\qquad \% 1\\
\mu_k & (\beta_k)^*\odot ((r^{(k)T}E_k x^{(k)}\oplus b_{k+1}),\\
x^{(k+1)} & \left(\begin{array}{l}
E_{k} y^{(k)}\mu_{k}\oplus x^{(k)}\\
\mu_{k}\end{array}\right).
\end{array}
\end{equation}

This expression yields the following generalized
version of the Levinson algorithm for solving linear symmetric 
Toeplitz systems \cite{GvL:99}:

\begin{algorithm}
\label{a:levi}
Bellman system with symmetric Toeplitz matrix.
\end{algorithm}
\noindent
${\bf function}\ y={\bf levinson}(r_0,r,b)\\
n={\bf size}(b)\\
y(1)=r_0^* \odot r(1);\quad x(1)=r_0^*\odot b(1)\\
\beta=r_0\qquad \% 1\\
\alpha=r_0^*\odot r(1)\qquad \% 1\\
{\bf for}\quad k=1:n-1\\
\-\beta=r_0\oplus r(1:k)\odot y(1:k)\qquad \% 2\\
\-\beta=\beta\oplus (\beta^*)^{-1}\odot\alpha^2\qquad \% 1\\
\-\mu=(r(k:-1:1)\odot x(1:k)\oplus b(k+1))\odot\beta^*\\
\- v(1:k)=x(1:k)\oplus\mu\odot y(k:-1:1)\\
\- x(1:k)=v(1:k)\\
\- x(k+1)=\mu\\
\- {\bf if}\quad k<n-1\\
\-\-\alpha=(r(k:-1:1)\odot y(1:k)\oplus r(k+1))\odot\beta^*\\
\-\- z(1:k)=y(1:k)\oplus\alpha\odot y(k:-1:1)\\
\-\- y(1:k)=z(1:k)\\
\-\- y(k+1)=\alpha\\
\-{\bf end}\\
{\bf end}$

The computational complexity of all 
methods described in this section is $O(n^2)$.

%\bibliographystyle{unsrt}
%\bibliography{toep1010}

\end{document}